\documentclass{amsart}
\usepackage{amssymb}
\usepackage[dvips]{graphicx}
\usepackage{graphics}
\usepackage{psfrag}
\begin{document}
\newtheorem{thm1}{Theorem}[section]
\newtheorem{lem1}[thm1]{Lemma}
\newtheorem{rem1}[thm1]{Remark}
\newtheorem{def1}[thm1]{Definition}
\newtheorem{cor1}[thm1]{Corollary}
\newtheorem{defn1}[thm1]{Definition}
\newtheorem{prop1}[thm1]{Proposition}
\newtheorem{ex1}[thm1]{Example}
\newtheorem{alg1}[thm1]{Algorithm}


\title[universal Gr\"{o}bner bases]{On the universal Gr\"{o}bner bases of toric ideals of graphs}
\author{Christos Tatakis}
\address{Department of Mathematics, University of Ioannina,
Ioannina 45110, Greece }
\email{chtataki@cc.uoi.gr}
\author{Apostolos Thoma }
\address{Department of Mathematics, University of Ioannina,
Ioannina 45110, Greece }
\email{athoma@uoi.gr}
\thanks{}


\date{}

\dedicatory{}

\begin{abstract}
\par
The universal Gr\"{o}bner  basis of $I$, is
a Gr\"{o}bner basis for $I$ with respect to all term orders simultaneously.
Let $I_G$ be the toric ideal of a graph $G$. We characterize in
graph theoretical terms the elements of the  universal Gr\"{o}bner basis of the toric ideal $I_G$.
 We provide a bound for the degree of the binomials in the universal Gr\"{o}bner  basis of the toric ideal of a graph.
Finally we give a family of examples of  circuits for which their true degrees are less than the degrees of some elements of the Graver basis.

\end{abstract}
\maketitle


\section{Introduction}

The universal Gr\"{o}bner basis of an ideal $I$
is the union of all reduced Gr\"{o}bner bases $G_<$ of the ideal $I$ as $<$ runs over all term orders.
The universal Gr\"{o}bner  basis  is a finite subset of $I$ and it is
a Gr\"{o}bner basis for $I$ with respect to all term orders simultaneously, see \cite{St}.
Universal Gr\"{o}bner bases exist for every ideal in $K[x_1,\dots , x_n]$. They were introduced
by V. Weispfenning \cite{W} and N. Schwartz \cite{S}.

\par  Let $A=\{{\bf a}_1,\ldots,{\bf a}_m\}\subseteq \mathbb{N}^n$
be a vector configuration in $\mathbb{Q}^n$ and
$\mathbb{N}A:=\{l_1{\bf a}_1+\cdots+l_m{\bf a}_m \ | \ l_i \in
\mathbb{N}\}$ the corresponding affine semigroup.  We grade the
polynomial ring $K[x_1,\ldots,x_m]$ over an arbitrary field $K$ by the
semigroup $\mathbb{N}A$ setting $\deg_{A}(x_i)={\bf a}_i$ for
$i=1,\ldots,m$. For ${\bf u}=(u_1,\ldots,u_m) \in \mathbb{N}^m$,
we define the $A$-{\em degree} of the monomial ${\bf x}^{{\bf
u}}:=x_1^{u_1} \cdots x_m^{u_m}$ to be \[ \deg_{A}({\bf x}^{{\bf
u}}):=u_1{\bf a}_1+\cdots+u_m{\bf a}_m \in \mathbb{N}A.\]  The
{\em toric ideal} $I_{A}$ associated to $A$ is the prime ideal
generated by all the binomials ${\bf x}^{{\bf u}}- {\bf x}^{{\bf
v}}$ such that $\deg_{A}({\bf x}^{{\bf u}})=\deg_{A}({\bf x}^{{\bf
v}})$, see \cite{St}. For such binomials, we set $\deg_A({\bf
x}^{{\bf u}}- {\bf x}^{{\bf v}}):=\deg_{A}({\bf x}^{{\bf u}})$.
An irreducible binomial ${\bf x}^{{\bf u}}- {\bf x}^{{\bf
v}}$ in $I_A$ is called {\em primitive} if there exists no other binomial
 ${\bf x}^{{\bf w}}- {\bf x}^{{\bf
z}}$ in $I_A$ such that ${\bf x}^{{\bf w}}$ divides $ {\bf x}^{{\bf
u}}$ and ${\bf x}^{{\bf z}}$ divides $ {\bf x}^{{\bf
v}}$.  The set of primitive binomials forms the Graver basis of $I_A$ and is denoted by
$Gr_A$.
An irreducible binomial is called {\em circuit} if it has minimal support. The set of circuits is
denoted by ${ C}_A$. The relation among the set of circuits,
 the Graver basis and the universal Gr\"{o}bner basis, which is denoted by ${ U}_A$, for a toric ideal $I_A$ is given by B. Sturmfels \cite{St}:

\begin{prop1} For any toric ideal $I_A$ we have ${ C}_A\subset { U}_A \subset  Gr_A$.

\end{prop1}

For toric ideals of graphs circuits were determined by R. Villarreal \cite[Proposition 4.2]{Vi}. The
 Graver basis of a toric ideal of a graph first have been studied by H. Ohsugi and T. Hibi \cite[Lemma 2.1]{OH} and the form of its elements
was determined by E. Reyes, Ch. Tatakis and A. Thoma \cite[Theorem 3.1 and Corollary 3.2]{RTT}.
While in \cite[Theorem 5.1]{LST} J. De Loera, B. Sturmfels and R. Thomas determined the universal Gr\"{o}bner basis for
toric ideals of graphs with less than nine vertices.
The purpose of this article is to determine the universal Gr\"{o}bner basis for the toric ideal of any graph.
In particular in section 2 we present some terminology, notations and results about the toric ideals of graphs. Section 3
contains the main result of the article which is a characterization of the binomials that belong to the universal  Gr\"{o}bner basis
of a toric ideal of a graph. Section 4 provides a degree bound for the binomials in the universal  Gr\"{o}bner basis of the toric ideal of a graph
and gives a family of examples of  circuits for which their true degrees are less than the degrees of some elements of the Graver basis.
Thus answering in the negative a conjecture by B. Sturmfels \cite[Conjecture 4.8]{St1}.

\section{Toric ideals of graphs}

\par
Let $G$ be a finite simple connected graph with vertices
$V(G)=\{v_{1},\ldots,v_{n}\}$ and edges $E(G)=\{e_{1},\ldots,e_{m}\}$.
A \emph{walk}  connecting $v_{1}\in V(G)$ and
$v_{q+1}\in V(G)$ is a finite sequence of the form
$$w=(\{v_{i_1},v_{i_2}\},\{v_{i_2},v_{i_3}\},\ldots,\{v_{i_q},v_{i_{q+1}}\})$$
with each $e_{i_j}=\{v_{i_j},v_{i_{j+1}}\}\in E(G)$.
\emph{Length}
of the walk $w$ is called the number $q$ of edges of the walk. An
even (respectively odd) walk is a walk of \emph{even} (respectively odd) length.
A walk
$w=(\{v_{i_1},v_{i_2}\},\{v_{i_2},v_{i_3}\},\ldots,\{v_{i_q},v_{i_{q+1}}\})$
is called \emph{closed} if $v_{i_{q+1}}=v_{i_1}$. A \emph{cycle}
is a closed walk
$$(\{v_{i_1},v_{i_2}\},\{v_{i_2},v_{i_3}\},\ldots,\{v_{i_q},v_{i_{1}}\})$$ with
$v_{i_k}\neq v_{i_j},$ for every $ 1\leq k < j \leq q$. Note that, although the graph $G$ has no multiple edges, the
same edge $e$ may appear more than once in a walk. In this case $e$ is
called {\em multiple edge of the walk $w$}.

Let $\mathbb{K}[e_{1},\ldots,e_{m}]$
 the polynomial ring in the $m$ variables $e_{1},\ldots,e_{m}$ over a field $\mathbb{K}$.  We
will associate each edge $e=\{v_{i},v_{j}\}\in E(G)$ with
$a_{e}=v_{i}+v_{j}$ in the free abelian group generated by the
vertices and let $A_{G}=\{a_{e}\ | \ e\in E(G)\}$. With $I_{G}$ we denote
 the toric ideal $I_{A_{G}}$ in
$\mathbb{K}[e_{1},\ldots,e_{m}]$.

Given an even closed walk of the graph $G$ $$w =(e_{i_1}, e_{i_2},\cdots,
e_{i_{2q}})$$ write
$$E^+(w)=\prod _{k=1}^{q} e_{i_{2k-1}},\ E^-(w)=\prod _{k=1}^{q} e_{i_{2k}}$$
and denote by $B_w$ the binomial
$$B_w=\prod _{k=1}^{q} e_{i_{2k-1}}-\prod _{k=1}^{q} e_{i_{2k}}.$$
It is easy to see that $B_w\in I_G$. Moreover, it is known that the toric ideal $I_G$
is generated by binomials of this form, see \cite{Vi}. For convenience
we denote by $\bf{w}$ the subgraph of $G$ with vertices the vertices of the
walk and edges the edges of the walk $w$. We call a walk
$w'=(e_{j_{1}},\dots,e_{j_{t}})$ a \emph{subwalk} of $w$ if
$e_{j_1}\cdots e_{j_t}| e_{i_1}\cdots e_{i_q}.$
An even
closed walk $w=(e_{i_1}, e_{i_2},\cdots, e_{i_{2q}})$ is said to
be primitive if there exists no even closed subwalk $\xi$ of $w$ of smaller
length such that $E^+(\xi)| E^+(w)$ and $E^-(\xi)| E^-(w)$. The walk $w$
is primitive if and only if the binomial $B_w$ is primitive.
Every even primitive walk $w=(e_{i_1},\ldots,e_{i_{2k}})$
partitions the set of edges in the two sets ${\bf w}^+= \{e_{i_j}|j \
{\it odd}\}, {\bf w}^-=\{e_{i_j}|j \ {\it even}\}$, otherwise the
binomial $B_w$ is not irreducible. While by $w^+$ we denote the exponent vector of the monomial $E^+(w)$
and by  $w^-$ the exponent vector of the monomial $E^-(w)$.

A {\em cut edge} (respectively {\em cut vertex}) is an edge (respectively vertex) of
the graph whose removal increases the number of connected
components of the remaining subgraph.  A graph is called {\em
biconnected} if it is connected and does not contain a cut
vertex. A {\em block} is a maximal biconnected subgraph of a given
graph $G$.
The edges of ${\bf w}^+$ are called odd edges of the walk and those of
${\bf w}^-$ even. Note that for a closed even walk whether an edge is even or
odd depends only on the edge that you start counting from. So it is
not important to identify whether an edge is even or odd but to separate the
edges in the two disjoint classes. A \emph{sink} of a block $B$ is a
common vertex of two odd or two even edges of the walk $w$ which
belong to the block $B$. In particular if $e$ is a cut edge of a
primitive walk then $e$ appears at least twice in the walk and
belongs either to ${\bf w}^+$ or ${\bf w}^-$. Therefore both vertices of $e$
are sinks. Sink is a property of the walk $w$ and not of the
underlying graph $\bf{w}$.

In the case of the toric ideals of graphs the following Theorems determine the form of the circuits and the
primitive binomials. R. Villarreal in
\cite[Proposition 4.2]{Vi} gave a necessary and sufficient
characterization of circuits:
\begin{thm1}\label{circuit}Let $G$ be a finite connected graph. The binomial $B\in I_{G}$ is circuit if and only if $B=B_{w}$ where
\begin{enumerate}
  \item $w$ is an even cycle or
  \item two odd cycles intersecting in exactly one vertex or
  \item two vertex disjoint odd cycles joined by a path.
\end{enumerate}
\end{thm1}

The next Theorem by  E. Reyes, Ch. Tatakis and A. Thoma describes  the form of the primitive binomials, i.e.
the elements $B_w\in I_G$ that belong to the Graver basis, \cite[Theorem 3.1]{RTT}.

\begin{thm1} \label{primitive}
Let $G$ a  graph and $w$ an even closed walk of $G$. The binomial $B_w$
is primitive if and only if
\begin{enumerate}
  \item every block of $\bf{w}$ is a cycle or a cut edge,
  \item every multiple edge of the walk $w$ is a double edge of the walk and a cut edge of $\bf{w}$,
  \item every cut vertex of $\bf{w}$ belongs to exactly two blocks and it is a sink of both.
\end{enumerate}
\end{thm1}
\begin{center}
\psfrag{D}{$B$}\psfrag{A}{Figure \ 1}
\includegraphics{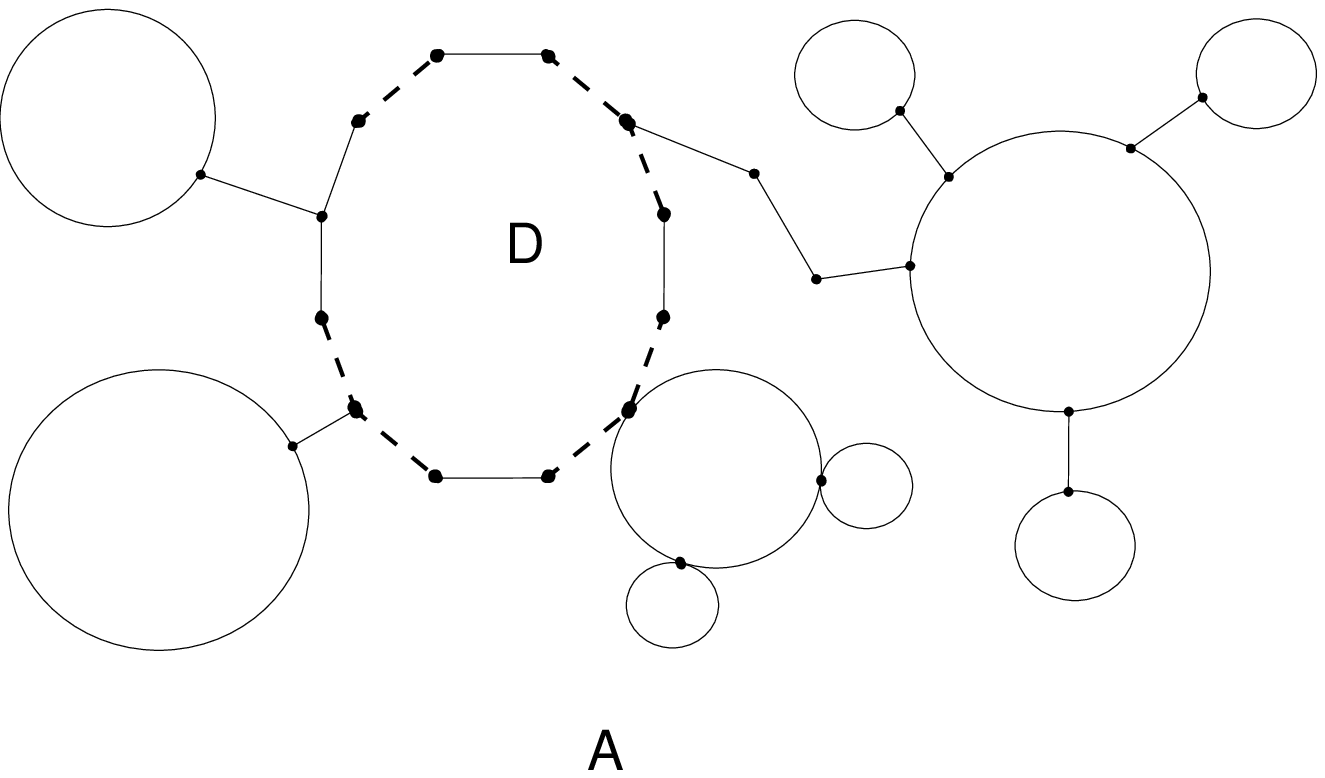}
\end{center}
Figure 1 shows a graph ${\bf w} $ of a primitive walk and a block B with four sinks.

\section{ Universal
Gr\"{o}bner bases }

In this section we will characterize the elements of the universal
Gr\"{o}bner basis of the toric ideal of a graph. The elements $B_w$ of the universal Gr\"{o}bner basis belong to the Graver basis, therefore their form
is determined by Theorem \ref{primitive}.  Let  $w=(e_{i_1}, e_{i_2},\cdots,
e_{i_{2q}})$ be a primitive walk then the blocks of the graph ${\bf w}$ are cyclic or they are cut edges. The simplest examble of a walk $w$ such that
$B_w$ is in the Graver basis but not in the universal Gr\"{o}bner basis is the one with degree 6 whose graph is in the figure 2.
\begin{center}
\psfrag{A}{Figure \ 2}
\includegraphics{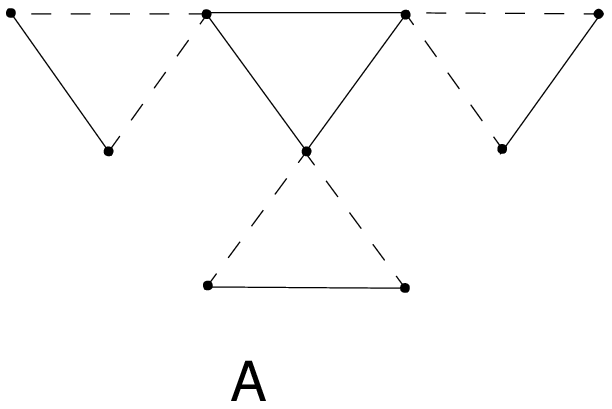}
\end{center}
The existence of this walk imply for $n\geq 9$
that ${ U}_{K_n} \not=  Gr_{K_n}$, where $K_n$ is the complete graph on $n$ vertices.
Note that in \cite{LST} J. De Loera, B. Sturmfels and R. Thomas  prove   that ${ C}_{K_n}= { U}_{K_n} =  Gr_{K_n}$ for $n\leq 7$ and
 ${ C}_{K_8}\not = { U}_{K_8} =  Gr_{K_8}$.  The reason for this walk not to be in the universal Gr\"{o}bner basis
is the existence of a {\em pure} cyclic block, the one in the center, that all of its edges are either in ${\bf w}^+$ or in  ${\bf w}^-$. In the next proposition
\ref{Prop} we will see that whenever a primitive walk $w$ has a block like that then the binomial $B_w$ is not in the universal Gr\"{o}bner basis.
In Theorem \ref{Theorem} we will see the converse, that is whenever an element $B_w$ is in the Graver basis but not in the universal Gr\"{o}bner basis
then $w$ has a pure cyclic block.

\begin{def1} A cyclic block $B$ of a primitive walk $w$ is called pure if all edges of $B$ are either in ${\bf w}^+$ or in  ${\bf w}^-$.
\end{def1}

\begin{prop1} \label{Prop} Let $w$ be an even primitive walk that has a pure cyclic block then $B_w$ does not belong to the universal
Gr\"{o}bner basis of $I_G$.
\end{prop1}
\textbf{Proof.} Suppose that $w$ has a pure cyclic block $B$ with edges $\epsilon_1,\dots, \epsilon_s$ which we can assume that belong to
${\bf w}^-$. Then the walk $w$ can be written in the
form $(w_1, \epsilon_1,\dots, w_s, \epsilon_s)$, where $w_i$ are subwalks of $w$ of odd length.
\begin{center}
\psfrag{A}{$w_{1}$}\psfrag{B}{$w_{2}$}\psfrag{C}{$w_{s}$}\psfrag{D}{$\epsilon_{1}$}\psfrag{E}{$\epsilon_{2}$}\psfrag{F}{$\epsilon_{s}$}\psfrag{H}{$B$}
\psfrag{G}{Figure \ 3}
\includegraphics{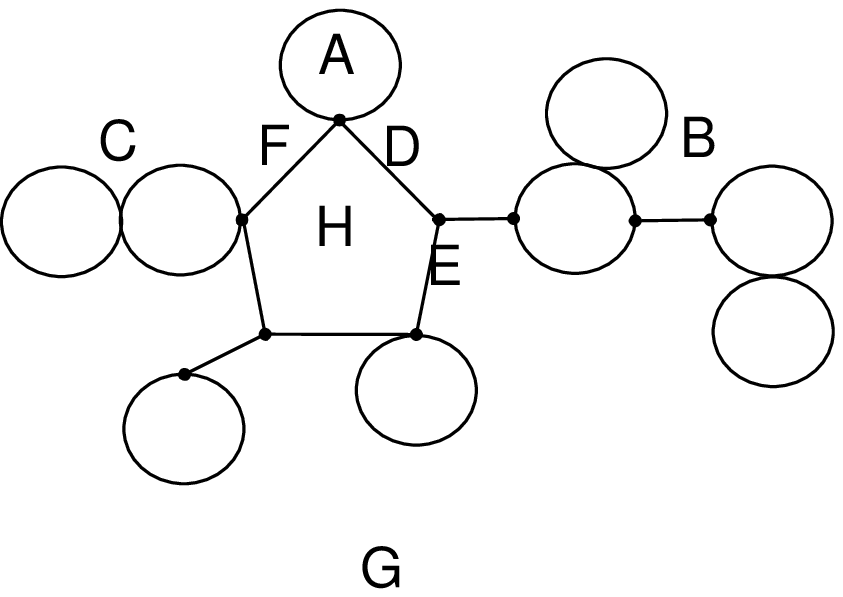}
\end{center}
For a subwalk $w_i $ we denote by
$$E^+(w_i)=\prod _{e_{i_{2k-1}}\in w_i} e_{i_{2k-1}},\ E^-(w_i)=\prod _{e_{i_{2k}}\in w_i} e_{i_{2k}}.$$
Then $B_w=E^+(w_1)E^+(w_2)\dots E^+(w_s)-\epsilon_1\epsilon_2 \dots \epsilon_s E^-(w_1)E^-(w_2)\dots E^-(w_s)$.
Look at the even walks $(w_i, \epsilon_i, w_{i+1}, \epsilon_{i})$ and the corresponding binomials
$F_i=E^+(w_i)E^+(w_{i+1})-\epsilon_i^2E^-(w_i)E^-(w_{i+1})\in I_G$, where $1\leq i\leq s-1$ and
$F_s=E^+(w_s)E^+(w_{1})-\epsilon_s^2E^-(w_s)E^-(w_{1})\in I_G$.\\
Suppose that $B_w$ belongs to a reduced Gr\"{o}bner basis for $I_G$ with respect to a term order $<$.
There are two cases.\\
First case: $E^+(w_1)E^+(w_2)\dots E^+(w_s)>\epsilon_1\epsilon_2 \dots \epsilon_s E^-(w_1)E^-(w_2)\dots E^-(w_s)$. Then necessarilly
 $E^+(w_i)E^+(w_{i+1})<\epsilon_i^2E^-(w_i)E^-(w_{i+1})$ for every $i$,
since $E^+(w_i)E^+(w_{i+1})$ divides $E^+(w_1)E^+(w_2)\dots E^+(w_s)$ and $F_i\in I_G$. \\
Multiplying all these inequalities for different i's we get  $$(E^+(w_1)E^+(w_2)\dots E^+(w_s))^2<(\epsilon_1\epsilon_2 \dots \epsilon_s E^-(w_1)E^-(w_2)\dots E^-(w_s))^2,$$
which is a contradiction.\\
Second case: $E^+(w_1)E^+(w_2)\dots E^+(w_s)<\epsilon_1\epsilon_2 \dots \epsilon_s E^-(w_1)E^-(w_2)\dots E^-(w_s)$.
In the case that $s=2k$ the binomial $G=\epsilon_1\epsilon_3\dots \epsilon_{2k-1}-\epsilon_2\epsilon_4\dots \epsilon_{2k}$
is in $I_G$ and both monomials of $G$
divide $\epsilon_1\epsilon_2 \dots \epsilon_{2k} E^-(w_1)E^-(w_2)\dots E^-(w_s)$,
a contradiction to the fact that $B_w$ belongs to the reduced Gr\"{o}bner basis.\\ In the case that $s=2k+1$ the binomials $G_i=E^+(w_i)\epsilon_{i+1}\epsilon_{i+3}\dots
\epsilon_{i+2k-1}-E^-(w_i)\epsilon_i\epsilon_{i+2}\dots \epsilon_{i+2k}$ are in $I_G$,
 where $\epsilon_j=\epsilon_l$ if $j\equiv l \ \textrm{mod} (2k+1)$. Therefore $E^+(w_i)\epsilon_{i+1}\epsilon_{i+3}\dots
\epsilon_{i+2k-1}>E^-(w_i)\epsilon_i\epsilon_{i+2}\dots \epsilon_{i+2k}$,
since $E^-(w_i)\epsilon_i\epsilon_{i+2}\dots \epsilon_{i+2k}$ divides  $\epsilon_1\epsilon_2 \dots
\epsilon_{2k} E^-(w_1)E^-(w_2)\dots E^-(w_s)$.
Multiplying them all and cancelling common factors
we get $$E^+(w_1)E^+(w_2)\dots E^+(w_s)>\epsilon_1\epsilon_2 \dots \epsilon_s E^-(w_1)E^-(w_2)\dots E^-(w_s),$$ a contradiction.
Therefore $B_w$ does not belong to any reduced Gr\"{o}bner basis of $I_G$ and thus also to the minimal universal
Gr\"{o}bner basis of $I_G$.
\hfill $\square$

\begin{def1} A primitive walk $w$ is called mixed if no cyclic block of $w$ is pure.
\end{def1}

The next Theorem is the main result of the article and describes the elements of the universal Gr\"{o}bner basis of $I_G$, for a general graph $G$. For any primitive walk
$w$ we construct a term order $<_w$ that
depends on $w$ to prove that a mixed primitive binomial belongs to the reduced Gr\"{o}bner basis with respect to this
 term order $<_w$. To prove it we will show that whenever one monomial of a binomial $B$ in $I_G$ divides one
of $E^+(w)$, $E^-(w)$ then the other monomial of $B$ is greater with respect to $<_w$ and does not divide either $E^+(w)$
or $E^-(w)$.

\begin{thm1} \label{Theorem} Let $w$ be a primitive walk. $B_w$ belongs to the universal Gr\"{o}bner basis of $I_G$ if and only if $w$ is mixed.

\end{thm1}
\textbf{Proof.} If $w$ is not mixed then it has a pure cyclic block and the result follows from Proposition \ref{Prop}. \\
Let $w$ be a mixed primitive walk.
We define a term order $<_w$ on $\mathbb{K}[e_1,\ldots,e_n]$, as an elimination order with the variables
that do not belong to $\bf{w}$ larger than the variables in ${\bf w}$. We order the first set of variables, with any term order and the second set of variables as follows:
Let $B_1, \dots B_{s_0}$ be any enumeration of all cyclic blocks of ${\bf w}$. Let $t_i^+$
 denotes the number of edges in ${\bf w}^+\cap B_i$ and
$t_i^-$ denotes the number of edges in ${\bf w}^-\cap B_i$. Let $W=(w_{ij})$ be the $(s_0)\times m$
 matrix where
$$
 w_{ij}= \left\{
  \begin{array}{ll}
    0, & \hbox{if $e_j \not\in B_i$,} \\
    t_i^-, & \hbox{if $e_j \in B_i \cap {\bf w}^+$,} \\
    t_i^+, & \hbox{if $e_j \in B_i \cap {\bf w}^-$}
  \end{array}
\right.
$$
and $m$ is the number of edges of ${\bf w}$.

Note that each column has at most one nonzero entry since each edge belongs to exactly one block of ${\bf w}$.
We say that $e^{ u}<_we^{v}$
if and only if the first nonzero coordinate of $W[{ u}-{ v}]$ is negative, otherwise, if $W[{ u}-{ v}]={\bf 0}$,
order them with any term order. Where $[u]$ is the vector $u$ written as a column vector. Note that for the walk $w$ we have  $W[w^+-w^-]={\bf 0}$.
Figure 4 shows a mixed primitive walk with their degrees $w_{ij}$.
\begin{center}
\psfrag{A}{2}\psfrag{B}{1}\psfrag{C}{1}\psfrag{D}{2}\psfrag{E}{2}\psfrag{F}{2}\psfrag{G}{2}\psfrag{H}{2}
\psfrag{I}{2}\psfrag{J}{2}\psfrag{K}{2}\psfrag{L}{4}\psfrag{M}{4}\psfrag{N}{1}\psfrag{O}{1}\psfrag{P}{2}
\psfrag{Q}{1}\psfrag{R}{1}\psfrag{S}{2}\psfrag{T}{0}\psfrag{U}{2}\psfrag{V}{3}\psfrag{W}{2}\psfrag{X}{3}
\psfrag{Y}{2}\psfrag{Z}{0}\psfrag{MN}{Figure \ 4}
\includegraphics{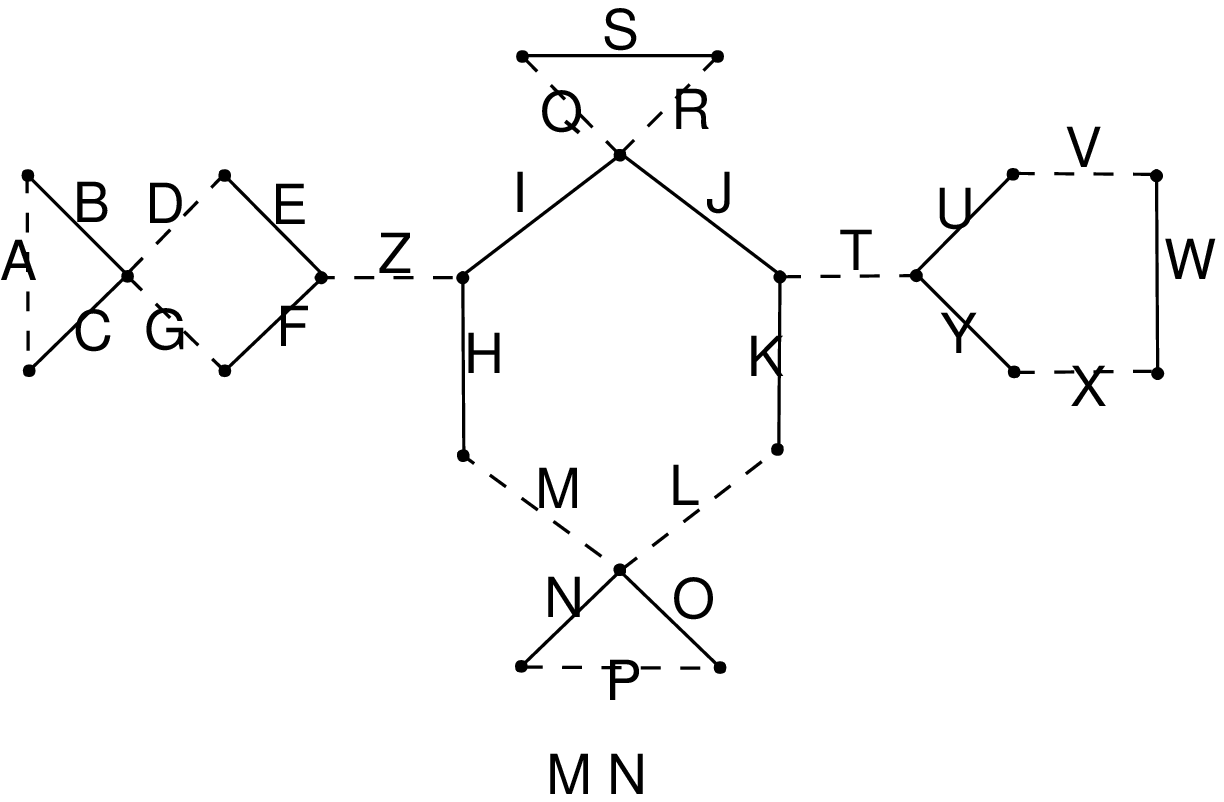}
\end{center}
We will prove that $B_w$ belongs to the reduced Gr\"{o}bner basis of $I_G$ with respect to the term order $<_w$.
It is enough to prove that whenever there exists a primitive binomial $B_z$ such that $E^+(z)|E^+(w)$ then $E^-(z)>_w E^+(z)$.
Note that $E^-(z)\nmid E^-(w)$ since $w$ is primitive and $E^-(z)\nmid E^+(w)$ since $w$ is mixed.
We remark that if ${\bf z}\nsubseteq{\bf w}$, since ${\bf z}^+\subset {\bf w}^+$,  there is an edge
 of ${\bf z}^-$ which is not an edge of ${\bf w}$.
But then $E^-(z)>_w E^+(z)$ since it is an elimination order. So we can suppose that ${\bf z}\subset {\bf w}$,
see also  \cite[proposition 4.13]{St}.

We claim that there exists at least one $i$, such that $B_i\cap{\bf z} \neq \emptyset$ and $B_i\cap {\bf z}^+
 \subsetneqq B_i\cap {\bf w}^+$.
Suppose not, then for every $i$, either $B_i\cap{\bf z} = \emptyset$ or $B_i\cap{\bf z}^+ =B_i\cap {\bf w}^+$
since $E^+(z)|E^+(w)$.
Let $B_i$ be a cyclic block such that $B_i\cap{\bf z}^+ =B_i\cap {\bf w}^+$, then $B_i\cap{\bf z}^- =B_i\cap {\bf w}^-$.
If not, then $B_i$ is not a block of ${\bf z}$ which
 implies that every edge $e$ in $B_i\cap {\bf z}^+$ is a cut edge of $\bf {z}$
and therefore $e$ is a double edge of $z$. But $B_i$ is a cyclic block
of ${\bf w}$ which means that every edge of $B_i$ is a single edge of $w$.
Therefore $e^2|E^+(z)$ and $e^2\nmid E^+(w)$ which is impossible since $E^+(z)|E^+(w)$.
Therefore $B_i\cap {\bf z} =B_i$ or $B_i\cap {\bf z} =\emptyset$.
This is obviously true also for blocks which are cut edges. But $z\not= w$, therefore  at least one block of ${\bf w}$ exists such that
$B_i\cap {\bf z} =\emptyset$ and at least one such that  $B_i\cap {\bf z} =B_i$.
The graph ${\bf w}$ is a graph of a walk so it is connected, so  two  adjacent blocks $B_j$ and $B_i$ exist
such that  $B_j\cap {\bf z} =\emptyset$ and   $B_i\cap {\bf z} =B_i$.
Let $v$ be the common cut vertex of $B_j$ and $B_i$. Then $2v$ appears in the degree of one of $deg_A(E^+(z))$, $deg_A(E^-(z))$ but not in the other one.
Therefore $B_{z}\not \in I_G$, a contradiction.

Let $i$ be the smallest integer such that $B_i\cap{\bf z} \neq \emptyset$ and $B_i\cap {\bf z}^+ \subsetneqq B_i\cap {\bf w}^+$.
Then according to the previous argument, the first $i-1$ coordinates of $W[z^+-z^-]$ are zero, since if
 $B_j\cap{\bf z} = \emptyset$ then  $w_j[z^+]=0=w_j[z^-]$ and  if $B_j\cap {\bf z}^+ = B_j\cap {\bf w}^+$,
then from the argument in the previous paragraph we have also $B_j\cap {\bf z}^- = B_j\cap {\bf w}^-$ and
then $w_j[z^+]=t_j^-t_j^+=w_j[z^-]$, where $w_j$ is the j-row of $W$.
For the block $B_i$ we have two cases: either $B_i\cap {\bf z} \neq B_i$ or $B_i\cap {\bf z} = B_i$.\\
First case: let $e\in B_i\cap {\bf z}$, then $e$ is a cut edge and then $e\in {\bf z}^-$,
otherwise $e$ is a double edge of $z$ and a simple of $w$, contradicting the fact that $E^+(z)|E^+(w)$.
So  every edge of $B_i\cap {\bf z}$ is in ${\bf z}^-$ and therefore  $w_i[z^+]=0$ and $w_i[z^-]>0$. Thus $E^-(z)>_wE^+(z)$.\\
Second case: $B_i\cap {\bf z}=B_i\Rightarrow B_i\cap {\bf z}^- = B_i\setminus (B_i\cap {\bf z}^+)$
 and since $B_i\cap {\bf z}^+ \subsetneqq B_i\cap {\bf w}^+$, we have  $w_i[z^+]<t_i^-t_i^+<w_i[z^-]$.
Therefore $E^-(z)>_wE^+(z)$.\\
We conclude that $B_w$ is in the reduced  Gr\"{o}bner basis with respect to the term order $<_w$ and thus it belongs to the universal
 Gr\"{o}bner basis of $I_A$.  \hfill $\square$

\section {Degree Bounds}

The number of elements in the universal Gr\"{o}bner basis is usually very large, for example in \cite{LST}
J. De Loera, B. Sturmfels and R. Thomas computed that the number
of the elements in the  universal Gr\"{o}bner basis of $I_{K_8}$ is 45570, where $K_n$ is the complete graph on $n$ vertices.
An estimate for the size of a  universal Gr\"{o}bner basis
can be a bound for the degrees of the elements in the  universal Gr\"{o}bner basis. Let $d_n$ be
the largest degree of a binomial in the universal Gr\"{o}bner basis for $I_{K_n}$
  In \cite{LST} J. De Loera, B. Sturmfels and R. Thomas proved that $d_n$ satisfies
$n-2\leq d_n\leq $
$\left(
 \begin{array}{c}
n\\
2\\
 \end{array}
\right).$ We will improve this result by proving that $d_n$ takes always the value $n-2$.

\begin{prop1} \label{bound} The largest degree $d_n$ of a binomial in the universal Gr\"{o}bner basis for $I_{K_n}$ is $d_n=n-2$,
for $n\geq 4$.
\end{prop1}
\textbf{Proof.} We will prove that the largest degree $d_n$ of a binomial in the Graver basis for $I_{K_n}$ is $d_n=n-2$ and it is attained
by a circuit, see also \cite{LST}. Circuits are always in the universal Gr\"{o}bner basis \cite{St}
therefore the result follows.
 Theorem 3.1 and Corollary 3.2 of \cite{RTT} imply that a primitive walk consists of blocks which are cut edges and
cyclic blocks, one if it is a cycle otherwise at least two. Let $w$ be a primitive walk and suppose that $w$ has $s_0$ cyclic blocks and $s_1$ cut edges.
Thus $s=s_0+s_1$ is the total number of blocks. From Theorem 3.1 of \cite{RTT} we know that there are exactly $s-1$ cut points and each one belongs
to exactly two blocks.
 Let $B_1, \dots, B_{s_0}$ be the cyclic blocks and
$t_i$ denotes the number of edges (vertices) of the cyclic block $B_i$.
Then the total number of vertices of ${\bf w}$ is
$$t_1+\dots +t_{s_0}+2s_1-(s-1)\leq n, $$ since the cut points are counted twice, see  Theorem 3.1 of \cite{RTT}. Two times the degree of $B_w$
is the sum of edges of the cyclic blocks $t_1+\dots +t_{s_0}$ plus two times the number of cut edges $s_1$, since
cut edges are double edges of the walk $w$ and edges of cycles are always single. Therefore
$$2deg(B_w)= t_1+\dots +t_{s_0}+2s_1\leq n+s-1.$$
So the largest degree is attained when the number of blocks of ${\bf w}$ is the largest possible
and if it is possible the walk $w$ pass through all the $n$ vertices, to have equality.
 But from $t_1+\dots +t_{s_0}+2s_1\leq n+s-1$ we get $s+(t_1-2)+\dots +(t_{s_0}-2)\leq n-1$.
 Note that
$(t_1-2)+\dots +(t_{s_0}-2)\geq 2$ since cyclic blocks have at least three vertices and the walk has at least two cyclic blocks, except
if $w$ is a cycle but in that case there is just one block.
Therefore $s\leq n-3$, but $s=n-3$ is possible with a circuit
with $n-5$ cut edges plus 2 cyclic blocks of three vertices each,  which has the maximal possible degree
$(n+(n-3)-1)/2=n-2.$ \hfill $\square$

Since any graph $G$ with $m$ vertices is a subgraph of the complete graph $K_m$ we have the following corollary.
\begin{cor1} \label{b} Let $G$ be a graph with $m$ vertices, $m\geq 4$.
The largest degree $d$ of a binomial in the universal Gr\"{o}bner basis for $I_{G}$ is $d\leq m-2$.
\end{cor1}

The knowledge of the form of the circuits \cite[Proposition 4.2]{Vi}, the elements of the Graver basis \cite[Theorem 3.1]{RTT},
the minimal systems of generators \cite[Theorem 4.13]{RTT}
 and the elements
of the universal Gr\"{o}bner basis of the toric ideal of a graph $G$, Theorem \ref{Theorem}
 and the variety and the easyness of description of graphs permit us   easily to produce examples of toric ideals having specific properties. For example
one can easily construct graphs such that the universal Gr\"{o}bner basis is equal to the Graver basis, just by avoiding creating pure blocks
in the elements of the Graver basis or making subdivisions in some  of the edges of pure blocks.  For other toric ideals that have this property see the recent work \cite{BHP} of T. Bogart, R. Hemmecke and S. Petrovi\'{c}.

In the proof of Proposition \ref{bound} the binomial that has the maximal degree in $I_{K_n}$ is a circuit.
  B. Sturmfels in his lecture at Santa Cruz (July 1995, see \cite{St1}),
made a conjecture that circuits have always the maximal degree among the elements of the Graver basis,
 but S. Hosten and R. Thomas gave a counterexample of a toric ideal such that the maximal degree of the elemements of the
Graver basis was 16 while the maximal degree of the circuits was 15, see \cite{St1}.
This example lead B. Sturmfels to alter the conjecture to: the degree of any element in the Graver basis $Gr_A$ of a toric ideal $I_A$
is bounded above by the maximal true degree of any circuit in ${ C}_A$, \cite[Conjecture 4.8]{St1}. Following \cite{St1} we define the true degree of a
circuit as follows:
Consider any circuit $C\in { C}_A$ and regard its support $supp(C)$ as a  subset of $A$. The lattice $\mathbb{Z}(supp(C))$ has finite index
in the lattice $R(supp(C))\cap \mathbb{Z}A$, which is called the index of the circuit $C$ and denoted by $ index(C)$. The {\em true degree} of the circuit $C$
is the product $degree(C)\cdot index(C)$.

Next we give a family of examples of  circuits for which their true degrees are less than the degrees of some elements of the Graver basis.
 Let us consider a graph $G$ consisting of a cycle of length $s$ and $s$ odd cycles of length $l$
each one attached to a vertex of
the initial cycle. Let $w$ be the walk that
pass from every edge of the graph $G$. The length of the walk $w$ is $ls+s=s(l+1)$, which is even.
\begin{center}
\psfrag{A}{$w_{1}$}\psfrag{B}{$w_{2}$}\psfrag{C}{$w_{s}$}\psfrag{D}{Figure \ 5}
\includegraphics{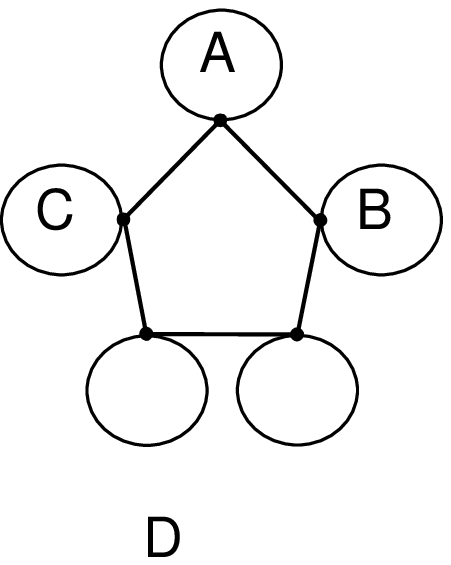}
\end{center}
Then $B_w$ is an element of the Graver basis of $I_G$, see
\cite{RTT}, and has degree $s(l+1)/2$. In the graph of $G$ there are a lot of circuits and the
 degree of the binomial corresponding to the longest one, which consists of two odd cycles joined by a path of length $s-1$, is
$(2l+2(s-1))/2=l+s-1$. Note that $s$, $l$, as lengths of cycles are greater than two, then $(s-2)(l-2)>0$ which
implies that $s(l+1)/2>l+s-1$. So there exists
an element $B_w$ in the Graver basis that has larger degree than any of the circuits, and it is easy to see that
the difference of the degrees can be made as
large as one wishes, by choosing big values for $l$ and $s$.
Note that an easy, but lengthy, computation of the the true degree of this circuit  shows that the true degree is
equal to the usual degree, therefore this example answers the question by B. Sturmfels \cite[Conjecture 4.8]{St1} in the negative.

Although the $B_w$ is in the Graver basis is not in the
universal Gr\"{o}bner basis, since it has a pure block, see Theorem   \ref{Theorem}. But if one takes a walk $w'$ such that ${\bf w}'$
consists of the
cycle in the center and any $(s-2)$
of the $s$ odd cycles then $w'$ is mixed and therefore $B_{w'}$ is in the universal Gr\"{o}bner basis and still the degree of $B_{w'}$
is bigger than any of
the degrees of circuits, for large $l$ and $s$.



\end{document}